\documentclass[12pt]{article}
\usepackage[a4paper,margin=1.3in]{geometry}
\usepackage{amssymb}
\usepackage{amsthm}
\usepackage{amsmath}
\usepackage{amsfonts}
\usepackage{graphicx}
\usepackage{float}
\usepackage[all]{xy}
\usepackage{scalerel}

\def\exp{\mathop{\mbox{\textsf{exp}}}\nolimits}

\newtheorem{theorem}{Theorem}
\newtheorem{lemma}{Lemma}
\newtheorem{definition}{Definition}
\newtheorem{proposition}{Proposition}
\newtheorem{corollary}{Corollary}
\newtheorem{example}{Example}
\newtheorem{remark}{Remark}

\begin{document}
\title{\textbf{Upper bound for flows of differential equations defined over normed rings}}
\author{Ronald Orozco}

\newcommand{\Addresses}{{

\textit{E-mail address}, R.~Orozco: \texttt{rj.orozco@uniandes.edu.co}
  
}}

\maketitle

\begin{abstract}
In this paper, we establish upper bounds for the flows of autonomous differential equations defined over complete normed rings. To achieve this, we compute explicit images of the autonomous operator acting on sequences in $\mathbb{R}^{\mathbb{N}}$ and employ a comparison principle based on Bell polynomials and the ring of convergent exponential generating functions. Finally, we illustrate our results with concrete applications to matrix differential equations and nonlinear initial value problems.
\end{abstract}
{\bf Keywords:} Autonomous operator, normed rings, dynamical systems, Bell polynomials, flows\\
{\bf Mathematics Subject Classification:} 47H99, 46H25, 37C10.

\section{Introduction}

The study of differential equations over abstract algebraic structures, such as integral domains and normed rings, provides a unified framework for analyzing systems ranging from classical real and complex equations to matrix and operator-valued dynamical systems. In \cite{orozco2}, the autonomous operator $\mathfrak{A}$ was introduced as an algebraic tool to construct formal power series solutions---or flows---of autonomous differential equations over integral domains via Hurwitz expansions.

While formal flows capture the algebraic structure of solutions, establishing their convergence and estimating their growth require analytic bounds. In this work, we extend the framework of \cite{orozco2} to complete normed rings. By combining the recursive structure of Bell polynomials with comparison series whose exact flows are known, we derive sharp upper bounds for the norm of the flow and determine explicit lower bounds for its radius of convergence.

This paper is organized as follows. In Section 2, we review the necessary preliminaries on Bell polynomials, Hurwitz expansions, the autonomous operator, and formal dynamical systems. In Section 3, we recall and compute explicit images of the autonomous operator for three fundamental sequences in $\mathbb{R}^{\mathbb{N}}$. In Section 4, we introduce the subring of convergent exponential generating functions over a complete normed ring and prove our main comparison theorem along with three explicit upper-bound corollaries. Finally, Section 5 presents concrete applications of our bounds to matrix differential equations and nonlinear systems. Throughout this paper, $\mathbb{N}$ denotes the monoid of non-negative integers.

\section{Preliminaries}

In this section, we review the definitions and basic properties of Bell polynomials, Hurwitz expansions, the autonomous operator, and dynamical systems over integral domains.

Bell polynomials are a fundamental combinatorial tool used to express the $n$-th derivative of a composition of two functions \cite{comtet}. Let $f$ and $g$ be two analytic functions with exponential power series representations $\sum_{n=0}^\infty a_n \frac{x^n}{n!}$ and $\sum_{n=0}^\infty b_n \frac{x^n}{n!}$, respectively, with $a_n, b_n \in \mathbb{C}$. Then
\begin{equation}
f(g(x)) = f(b_0) + \sum_{n=1}^{\infty} Y_n(b_1, \dots, b_n; a_1, \dots, a_n) \frac{x^n}{n!},
\end{equation}
where $Y_n$ is the $n$-th exponential Bell polynomial. For instance, the first three polynomials are
\begin{align*}
Y_1(b_1; a_1) &= a_1 b_1, \\
Y_2(b_1, b_2; a_1, a_2) &= a_1 b_2 + a_2 b_1^2, \\
Y_3(b_1, b_2, b_3; a_1, a_2, a_3) &= a_1 b_3 + 3 a_2 b_1 b_2 + a_3 b_1^3.
\end{align*}
Explicitly, Bell polynomials are given by the Fa\`a di Bruno formula \cite{comtet}:
\begin{equation}
Y_n(b_1, \dots, b_n; a_1, \dots, a_n) = \sum_{k=1}^{n} B_{n,k}(b_1, \dots, b_{n-k+1}) a_k,
\end{equation}
where
\begin{equation}
B_{n,k} = \sum_{\substack{j_1 + 2j_2 + \dots + n j_n = n \\ j_1 + j_2 + \dots + j_n = k}} \frac{n!}{j_1! j_2! \cdots j_n!} \left(\frac{b_1}{1!}\right)^{j_1} \left(\frac{b_2}{2!}\right)^{j_2} \cdots \left(\frac{b_n}{n!}\right)^{j_n},
\end{equation}
with the sum taken over all partitions of $n$ into $k$ parts, where $j_h$ denotes the number of parts of size $h$. Note that every coefficient in $B_{n,k}$ is non-negative.

Let $R$ be a commutative integral domain with identity, and let $\mathrm{Exp}_R[x]$ denote the ring of formal exponential power series of the form $\sum_{n=0}^\infty a_n \frac{x^n}{n!}$ with $a_n \in R$. Let $\delta : \mathrm{Exp}_R[x] \to \mathrm{Exp}_R[x]$ denote the formal derivative defined by $\delta\left(\sum_{n=0}^\infty a_n \frac{x^n}{n!}\right) = \sum_{n=0}^\infty a_{n+1} \frac{x^n}{n!}$, which satisfies $\delta(f+g) = \delta f + \delta g$ and $\delta(fg) = f\delta g + (\delta f)g$. We write $\delta^k f = \delta(\delta^{k-1} f)$ with $\delta^0 f = f$.

\begin{definition}
The \emph{Hurwitz expansion} of $f(x) \in \mathrm{Exp}_R[x]$, denoted $\Delta f(x)$, is the sequence of derivatives $(\delta^n f(x))_{n \in \mathbb{N}}$. Let $\Delta \mathrm{Exp}_R[x]$ denote the ring of all Hurwitz expansions equipped with componentwise addition and the Hurwitz convolution product $*$:
\begin{enumerate}
    \item $\Delta f + \Delta g = \Delta(f + g)$,
    \item $\Delta f * \Delta g = \Delta(fg)$, where the $n$-th component is given by
    \[
    (\Delta f(x) * \Delta g(x))_n = \sum_{k=0}^{n} \binom{n}{k} \delta^k f(x) \, \delta^{n-k} g(x).
    \]
\end{enumerate}
\end{definition}

\begin{definition}[\cite{orozco2}]
Let $\mathrm{Exp}_R[x]^{\mathbb{N}}$ be the algebra of sequences over $\mathrm{Exp}_R[x]$ with componentwise operations. The \emph{autonomous operator} $\mathfrak{A} : \Delta \mathrm{Exp}_R[x] \to \mathrm{Exp}_R[x]^{\mathbb{N}}$ is the nonlinear map defined by $\mathfrak{A}(\Delta f) = \left(A_n([\delta^{n-1} f])\right)_{n \geq 1}$, where $[\delta^n f] = (f, \delta f, \dots, \delta^n f)$ and the polynomials $A_n$ are defined recursively by:
\begin{align}
A_1([\delta^0 f]) &= f, \\
A_{n+1}([\delta^n f]) &= Y_n\left(A_1([\delta^0 f]), A_2([\delta^1 f]), \dots, A_n([\delta^{n-1} f]); \delta^1 f, \dots, \delta^n f\right), \quad n \geq 1.
\end{align}
\end{definition}

\begin{definition}[\cite{orozco2}]
For $\Delta f(x) \in \Delta \mathrm{Exp}_R[x]$, the formal \emph{dynamical system} $\Phi : R \times R \times \Delta \mathrm{Exp}_R[x] \to \mathrm{Exp}_S[t]$ (with $S = \mathrm{Exp}_R[x]^{\mathbb{N}}$) associated with $\Delta f(x)$ is defined by
\begin{equation}
\Phi(t, x, \Delta f(x)) = x + \sum_{n=1}^{\infty} A_n([\delta^{n-1} f(x)]) \frac{t^n}{n!}.
\end{equation}
For any ring extension $R'$ of $R$ and any $c \in R'$ where the evaluation map $\epsilon_c : \mathrm{Exp}_R[x] \to R'$ is well-defined, we extend $\epsilon_c$ to $\Phi$ by setting
\begin{equation}
\epsilon_c \Phi(t, x, \Delta f(x)) = c + \sum_{n=1}^{\infty} \epsilon_c\left(A_n([\delta^{n-1} f(x)])\right) \frac{t^n}{n!}.
\end{equation}
The trajectory of $\Phi$ starting at $c$ is given by $\Gamma_c := \{ \epsilon_c \Phi(t, x, \Delta f(x)) \}$.
\end{definition}

\section{Some explicit images of $\mathfrak{A}$}

In \cite{orozco1}, certain images of the autonomous operator $\mathfrak{A}$ were studied to obtain trajectories $\Gamma_0$ for the autonomous flow $\delta_t \Phi = f(\Phi)$. Below, we summarize and derive the images of three standard sequences in $\mathbb{R}^{\mathbb{N}}$.

\begin{proposition}\label{prop:images}
Let $a \in \mathbb{R}$ with $a \neq 0, 1$, and denote the falling factorial by $(a)_0 = 1$ and $(a)_k = a(a-1)\cdots(a-k+1)$ for $k \geq 1$. Then:
\begin{align}
\mathfrak{A}\left((a^n)_{n \geq 0}\right) &= \left(a^{n-1}(n-1)!\right)_{n \geq 1}, \label{eq:img1} \\
\mathfrak{A}\left((n!)_{n \geq 0}\right) &= \left((-1)^{n+1} 2^n \binom{1/2}{n} n!\right)_{n \geq 1}, \label{eq:img2} \\
\mathfrak{A}\left(((a)_n)_{n \geq 0}\right) &= \left(\binom{\frac{1}{1-a}}{n} (1-a)^n n!\right)_{n \geq 1}. \label{eq:img3}
\end{align}
\end{proposition}

\begin{proof}
Consider the autonomous initial value problem $y'(t) = f(y(t))$ with $y(0) = 0$, where $f$ is analytic in a neighborhood of $0$ with $f(0) = 1$. By the recursive definition of $A_n$ and Fa\`a di Bruno's formula, the unique analytic solution is given by the exponential generating function
\[
y(t) = \sum_{n=1}^{\infty} A_n\left(f(0), f'(0), \dots, f^{(n-1)}(0)\right) \frac{t^n}{n!}.
\]
Separating variables yields $\int_0^{y(t)} \frac{du}{f(u)} = t$. We evaluate this integral for three choices of $f$:
\begin{enumerate}
    \item For $f(z) = e^{az}$, we have $f^{(i)}(0) = a^i$. Solving $\int_0^y e^{-au} \, du = t$ gives
    \[
    y(t) = \frac{1}{a} \ln\left(\frac{1}{1-at}\right) = \sum_{n=1}^{\infty} a^{n-1} \frac{t^n}{n} = \sum_{n=1}^{\infty} a^{n-1}(n-1)! \frac{t^n}{n!}, \quad |at| < 1.
    \]
    \item For $f(z) = \frac{1}{1-z}$, we have $f^{(i)}(0) = i!$. Solving $\int_0^y (1-u) \, du = t$ gives $y - \frac{y^2}{2} = t$, whose solution satisfying $y(0)=0$ is
    \[
    y(t) = 1 - \sqrt{1-2t} = \sum_{n=1}^{\infty} (-1)^{n+1} 2^n \binom{1/2}{n} t^n, \quad |2t| < 1.
    \]
    \item For $f(z) = (1+z)^a$, we have $f^{(i)}(0) = (a)_i$. Solving $\int_0^y (1+u)^{-a} \, du = t$ gives $\frac{(1+y)^{1-a}-1}{1-a} = t$, hence
    \[
    y(t) = [1 + (1-a)t]^{\frac{1}{1-a}} - 1 = \sum_{n=1}^{\infty} \binom{\frac{1}{1-a}}{n} (1-a)^n t^n, \quad |(1-a)t| < 1.
    \]
\end{enumerate}
Equating coefficients of $\frac{t^n}{n!}$ in each expansion yields \eqref{eq:img1}, \eqref{eq:img2}, and \eqref{eq:img3}.
\end{proof}

\begin{corollary}\label{cor:flows_real}

For the functions $f(x) = e^{\alpha x}$, $f(x) = \frac{1}{1-x}$, and $f(x) = (1+x)^\alpha$ ($\alpha \neq 0, 1$), the flows evaluated at $0$ are:
\begin{align*}
\epsilon_0 \Phi\left(t, x, \Delta(e^{\alpha x})\right) &= \frac{1}{\alpha} \ln\left(\frac{1}{1-\alpha t}\right), \quad |\alpha t| < 1, \\
\epsilon_0 \Phi\left(t, x, \Delta\left(\frac{1}{1-x}\right)\right) &= 1 - \sqrt{1-2t}, \quad |2t| < 1, \\
\epsilon_0 \Phi\left(t, x, \Delta((1+x)^\alpha)\right) &= [1 + (1-\alpha)t]^{\frac{1}{1-\alpha}} - 1, \quad |(1-\alpha)t| < 1.
\end{align*}
\end{corollary}

\section{Upper bounds for $\Phi$ over normed rings}

We now use the comparison results from Section 3 to establish upper bounds for the flow of a one-dimensional dynamical system defined over a complete normed ring. Throughout this section, $R$ denotes a commutative normed integral domain with identity.
\begin{definition}[\cite{gelfand}]
A \emph{normed ring} is a ring $R$ endowed with a real-valued function $|\cdot| : R \to \mathbb{R}_{\geq 0}$, called a norm, satisfying for all $x, y \in R$:
\begin{enumerate}
\item $|x| = 0$ if and only if $x = 0$,
\item $|x + y| \leq |x| + |y|$,
\item $|-x| = |x|$,
\item $|xy| \leq |x| |y|$.
\end{enumerate}
A normed ring $R$ is \emph{complete} if every Cauchy sequence in $R$ converges to an element of $R$.
\end{definition}
Let $R$ be a complete normed ring. Consider the subset of $\mathrm{Exp}_R[x]$ defined by
\(\mathrm{Exp}_{R}\{x\}=\left\{\sum _{n=0}^{\infty }a_{n}\frac{x^{n}}{n!}:a_{n}\in R,\;\lim _{n\rightarrow \infty }\left\|\frac{a_{n}}{n!}\right\|=0\right\}.\)
Since the sequence $\left(\left|\frac{a_n}{n!}\right|\right)_{n \in \mathbb{N}}$ converges to $0$, it is bounded and attains its maximum. Thus, we can endow $\mathrm{Exp}_R{x}$ with the supremum norm
\(\|f\|=\max _{n\in \mathbb{N}}\left\|\frac{a_{n}}{n!}\right\|.\)
To ensure stability under arbitrary differentiation $\delta^k$, let $\mathrm{Exp}R^{\mathrm{an}}{x} \subset \mathrm{Exp}R{x}$ denote the subring of series $f(x) = \sum{n=0}^\infty a_n \frac{x^n}{n!}$ such that $\sum{n=0}^\infty \left|\frac{a_n}{n!}\right| < \infty$ and $\lim_{n\to\infty} \frac{(n+k)!}{n!} \left|\frac{a_{n+k}}{(n+k)!}\right| = 0$ for all $k \geq 0$. We denote its Hurwitz expansions by $\Delta \mathrm{Exp}_R^{\mathrm{an}}{x}$.
\begin{lemma}\label{lem:subring}
$\Delta \mathrm{Exp}_R^{\mathrm{an}}{x}$ is a subring of $\Delta \mathrm{Exp}_R[x]$ containing $\Delta R$.
\end{lemma}
\begin{proof}
Let $f(x) = \sum_{n=0}^\infty a_n \frac{x^n}{n!}$ and $g(x) = \sum_{n=0}^\infty b_n \frac{x^n}{n!}$ belong to $\mathrm{Exp}R^{\mathrm{an}}{x}$. Set $\alpha_n = a_n / n!$ and $\beta_n = b_n / n!$, so that $\sum{n=0}^\infty |\alpha_n| < \infty$ and $\sum_{n=0}^\infty |\beta_n| < \infty$. Clearly, $\sum_{n=0}^\infty |\alpha_n + \beta_n| < \infty$, so $f + g \in \mathrm{Exp}_R^{\mathrm{an}}{x}$.
For the product $fg = \sum_{n=0}^\infty c_n \frac{x^n}{n!}$, observe that
\(\gamma _{n}:=\frac{c_{n}}{n!}=\frac{1}{n!}\sum _{k=0}^{n}{n \choose k}a_{k}b_{n-k}=\sum _{k=0}^{n}\alpha _{k}\beta _{n-k}.\)
Since $\sum |\alpha_n| < \infty$ and $\sum |\beta_n| < \infty$, the Cauchy product $\sum_{n=0}^\infty |\gamma_n|$ converges by Mertens' theorem, which implies $\lim_{n\to\infty} |\gamma_n| = 0$. By the Leibniz rule, $\delta^k(fg) = \sum_{j=0}^k \binom{k}{j} (\delta^j f)(\delta^{k-j} g)$, so the same summability property holds for every derivative $\delta^k(fg)$. Hence $\Delta \mathrm{Exp}_R^{\mathrm{an}}{x}$ is closed under addition and Hurwitz convolution.
\end{proof}
\begin{remark}
For any complete normed ring extension $R'$ of $R$ and any $c \in R'$ with $|c| \leq 1$, the evaluation map $\epsilon_c : \mathrm{Exp}_R^{\mathrm{an}}{x} \to R'$ given by
\(\epsilon _{c}(f(x))=f(c)=\sum _{n=0}^{\infty }a_{n}\frac{c^{n}}{n!}\)
is a well-defined, continuous ring homomorphism because $\sum_{n=0}^\infty \left|a_n \frac{c^n}{n!}\right| \leq \sum_{n=0}^\infty \left|\frac{a_n}{n!}\right| |c|^n < \infty$.
\end{remark}

\begin{theorem}[Comparison Principle]\label{thm:comparison}
Let $\Delta f(x) \in \Delta \mathrm{Exp}R^{\mathrm{an}}{x}$ and let $g(t) = \sum{n=0}^\infty \alpha_n \frac{t^n}{n!} \in \mathrm{Exp}_{\mathbb{R}}{t}$ be a real comparison series with $\alpha_n \geq 0$ for all $n \geq 0$. Suppose there exists $c \in R'$ (with $|c| \leq 1$) such that
\(\|\epsilon _{c}(\delta ^{n}f(x))\|\le \alpha _{n}\quad \text{for\ all\ }n\ge 0,\)
and suppose that $\Phi(t, 0, \Delta g(0))$ converges for $0 \leq t < T$. Then, for all $n \geq 1$,
\begin{equation}\label{eq:comp_coeff}
|\epsilon_c(A_n([\delta^{n-1} f(x)]))| \leq A_n([\alpha_{n-1}]),
\end{equation}
where $[\alpha_{n-1}] = (\alpha_0, \alpha_1, \dots, \alpha_{n-1})$, and for all $t \in \mathbb{R}$ with $|t| < T$,
\begin{equation}\label{eq:comp_flow}
|\epsilon_c \Phi(t, x, \Delta f(x)) - c| \leq \Phi(|t|, 0, \Delta g(0)).
\end{equation}
\end{theorem}
\begin{proof}
We prove \eqref{eq:comp_coeff} by induction on $n$. For $n = 1$, we have
\(\|\epsilon _{c}(A_{1}([\delta ^{0}f(x)]))\|=\|\epsilon _{c}(f(x))\|\le \alpha _{0}=A_{1}([\alpha _{0}]).\)
Assume the inequality holds for all positive integers $1 \leq m \leq n$. For $n + 1$, using the submultiplicativity of the norm, the non-negativity of the coefficients $B_{n,k}$, and the induction hypothesis, we obtain
\begin{align*}
&|\epsilon_c(A_{n+1}([\delta^n f(x)]))|\\
&= \left| \sum_{k=1}^{n} \sum_{\substack{j_1 + \dots + n j_n = n \ j_1 + \dots + j_n = k}} \frac{n!}{j_1! \cdots j_n!} \left[\frac{\epsilon_c(A_1)}{1!}\right]^{j_1} \cdots \left[\frac{\epsilon_c(A_n)}{n!}\right]^{j_n} \epsilon_c(\delta^k f) \right| \\
&\leq \sum_{k=1}^{n} \sum_{\substack{j_1 + \dots + n j_n = n \ j_1 + \dots + j_n = k}} \frac{n!}{j_1! \cdots j_n!} \left[\frac{|\epsilon_c(A_1)|}{1!}\right]^{j_1} \cdots \left[\frac{|\epsilon_c(A_n)|}{n!}\right]^{j_n} |\epsilon_c(\delta^k f)| \\
&\leq \sum_{k=1}^{n} \sum_{\substack{j_1 + \dots + n j_n = n \ j_1 + \dots + j_n = k}} \frac{n!}{j_1! \cdots j_n!} \left[\frac{A_1([\alpha_0])}{1!}\right]^{j_1} \cdots \left[\frac{A_n([\alpha_{n-1}])}{n!}\right]^{j_n} \alpha_k \\
&= A_{n+1}([\alpha_n]).
\end{align*}
This completes the induction. Now, for any $t \in \mathbb{R}$ with $|t| < T$, since $A_n([\alpha_{n-1}]) \geq 0$, we have
\begin{align*}
|\epsilon_c \Phi(t, x, \Delta f(x)) - c| &= \left| \sum_{n=1}^{\infty} \epsilon_c(A_n([\delta^{n-1} f(x)])) \frac{t^n}{n!} \right| \\
&\leq \sum_{n=1}^{\infty} |\epsilon_c(A_n([\delta^{n-1} f(x)]))| \frac{|t|^n}{n!} \\
&\leq \sum_{n=1}^{\infty} A_n([\alpha_{n-1}]) \frac{|t|^n}{n!} = \Phi(|t|, 0, \Delta g(0)).
\end{align*}
In particular, the series defining $\epsilon_c \Phi(t, x, \Delta f(x))$ converges absolutely in the complete normed ring $R'$ for all $|t| < T$.
\end{proof}
Combining Theorem \ref{thm:comparison} with Proposition \ref{prop:images} and Corollary \ref{cor:flows_real}, we immediately obtain three explicit upper bounds.
\begin{corollary}\label{cor:exp}
Let $\Delta f(x) \in \Delta \mathrm{Exp}_R^{\mathrm{an}}{x}$ be such that $|\epsilon_c(\delta^n f(x))| \leq a^n$ for all $n \geq 0$, for some $c \in R'$ ($|c| \leq 1$) and $a > 0$. Then, for all $n \geq 1$,
\(\|\epsilon _{c}(A_{n}([\delta ^{n-1}f(x)]))\|\le a^{n-1}(n-1)!,\)
and for all $t \in \mathbb{R}$ with $a|t| < 1$,
\(\|\epsilon _{c}\Phi (t,x,\Delta f(x))-c\|\le \frac{1}{a}\ln \left(\frac{1}{1-a|t|}\right).\)
\end{corollary}
\begin{corollary}\label{cor:geom}
Let $\Delta f(x) \in \Delta \mathrm{Exp}_R^{\mathrm{an}}{x}$ be such that $|\epsilon_c(\delta^n f(x))| \leq n!$ for all $n \geq 0$, for some $c \in R'$ ($|c| \leq 1$). Then, for all $n \geq 1$,
\(\|\epsilon _{c}(A_{n}([\delta ^{n-1}f(x)]))\|\le (-1)^{n+1}2^{n}{1/2 \choose n}n!,\)
and for all $t \in \mathbb{R}$ with $2|t| < 1$,
\(\|\epsilon _{c}\Phi (t,x,\Delta f(x))-c\|\le 1-\sqrt{1-2|t|}.\)
\end{corollary}
\begin{corollary}\label{cor:bin}
Let $\Delta f(x) \in \Delta \mathrm{Exp}_R^{\mathrm{an}}{x}$ and let $c \in R'$ with $|c| \leq 1$.
\begin{enumerate}
\item \emph{(Polynomial comparison)} If $m \geq 2$ is an integer and $|\epsilon_c(\delta^n f(x))| \leq (m)_n = m(m-1)\cdots(m-n+1)$ for all $n \geq 0$, then for $(m-1)|t| < 1$:
\[\|\epsilon _{c}\Phi (t,x,\Delta f(x))-c\|\le [1-(m-1)|t|]^{-\frac{1}{m-1}}-1.\]
\item \emph{(General rising-factorial comparison)} For any real $a > 0$, if $|\epsilon_c(\delta^n f(x))| \leq a(a+1)\cdots(a+n-1)$ for all $n \geq 0$, then for $(a+1)|t| < 1$:
\[\|\epsilon _{c}\Phi (t,x,\Delta f(x))-c\|\le 1-[1-(a+1)|t|]^{\frac{1}{a+1}}.\]
\end{enumerate}
\end{corollary}

\section{Applications}

In this section, we illustrate how Theorem \ref{thm:comparison} and Corollaries \ref{cor:exp}--\ref{cor:bin} yield explicit convergence intervals and growth estimates for nonlinear differential equations over complete normed rings and Banach algebras.
\begin{example}[Matrix exponential differential equation]
Let $R' = \mathcal{M}_k(\mathbb{C})$ be the complete normed ring of $k \times k$ complex matrices equipped with a submultiplicative matrix norm $|\cdot|$ (such as the spectral norm or Frobenius norm, normalized so that $|I| = 1$), and let $R = \mathbb{C} I \cong \mathbb{C}$. Consider the matrix initial value problem
\begin{equation}\label{eq:mat_exp}
Y'(t) = B \exp(A Y(t)), \qquad Y(0) = C,
\end{equation}
where $A, B, C \in \mathcal{M}_k(\mathbb{C})$ commute pairwise, $|B| \leq 1$, $|C| \leq 1$, and $|A| \leq a$ for some $a > 0$. Here, the flow corresponds to evaluating the Hurwitz expansion of $f(x) = B e^{Ax}$ at $c = C$. Since $A$ and $B$ commute, the $n$-th derivative of $f$ evaluated at $C$ is
\(\epsilon _{C}(\delta ^{n}f(x))=BA^{n}\exp (AC).\)
Using the submultiplicativity of the norm, we estimate
\(\|\epsilon _{C}(\delta ^{n}f(x))\|\le \|B\|\|A\|^{n}\|\exp (AC)\|\le e^{a\|C\|}a^{n}.\)
Rescaling time by $\tau = e^{a|C|} t$, the hypotheses of Corollary \ref{cor:exp} are satisfied. Consequently, the unique analytic solution $Y(t)$ of \eqref{eq:mat_exp} exists on the interval $|t| < \frac{1}{a e^{a|C|}}$ and satisfies the explicit norm bound
\(\|Y(t)-C\|\le \frac{1}{a}\ln \left(\frac{1}{1-ae^{a\|C\|}|t|}\right).\)
\end{example}
\begin{example}[Matrix Riccati-type rational flow]
Let $R' = \mathcal{M}_k(\mathbb{C})$ as above, and consider the nonlinear matrix differential equation
\begin{equation}\label{eq:mat_riccati}
Y'(t) = (I - M Y(t))^{-1}, \qquad Y(0) = C,
\end{equation}
where $M, C \in \mathcal{M}_k(\mathbb{C})$ commute, $|C| \leq 1$, and $|M| \leq \mu$ with $\mu |C| < 1$. The vector field is $f(x) = (I - Mx)^{-1}$, whose $n$-th derivative at $C$ is
\(\epsilon _{C}(\delta ^{n}f(x))=n!\,M^{n}(I-MC)^{-(n+1)}.\)
Since $|(I - MC)^{-1}| \leq (1 - \mu |C|)^{-1}$, setting $\kappa = (1 - \mu|C|)^{-1}$ and $\lambda = \mu \kappa$ gives
\(\|\epsilon _{C}(\delta ^{n}f(x))\|\le \kappa \lambda ^{n}n!\qquad \text{for\ all\ }n\ge 0.\)
Applying Corollary \ref{cor:geom} after rescaling $t \mapsto \kappa \lambda t$ (noting that scaling $f$ by $\kappa$ and $x$ by $\lambda$ scales $A_n$ by $\kappa^n \lambda^{n-1}$), we find that the solution $Y(t)$ of \eqref{eq:mat_riccati} converges for $|t| < \frac{(1 - \mu|C|)^2}{2\mu}$ and obeys
\(\|Y(t)-C\|\le \frac{1}{\lambda }\left(1-\sqrt{1-2\kappa \lambda |t|}\right)=\frac{1-\mu \|C\|}{\mu }\left(1-\sqrt{1-\frac{2\mu |t|}{(1-\mu \|C\|)^{2}}}\right).\)
Notice that when $k=1$, $M = \mu > 0$, and $C = 0$, this bound reduces to $\frac{1}{\mu}(1 - \sqrt{1 - 2\mu|t|})$, which is the exact closed-form solution of $y'(t) = (1 - \mu y(t))^{-1}$ with $y(0) = 0$. This confirms that the upper bound in Corollary \ref{cor:geom} is sharp.
\end{example}
\begin{example}[Polynomial flows in commutative Banach algebras]
Let $\mathcal{B}$ be a unital commutative Banach algebra (for instance, $C([0,1])$ with the supremum norm or the Wiener algebra $A(\mathbb{T})$ of absolutely convergent Fourier series). Consider the polynomial differential equation
\(u^{\prime }(t)=(1+u(t))^{m},\qquad u(0)=u_{0}\in \mathcal{B},\)
where $m \geq 2$ is an integer and $|1 + u_0| \leq 1$. Here $\epsilon_{u_0}(\delta^n f(x)) = (m)n (1 + u_0)^{m-n}$ for $0 \leq n \leq m$, and $0$ for $n > m$. Thus, $|\epsilon{u_0}(\delta^n f(x))| \leq (m)_n$ for all $n \geq 0$. By Corollary \ref{cor:bin}(1), the flow $u(t)$ exists for all $|t| < \frac{1}{m-1}$ and satisfies
\(\|u(t)-u_{0}\|\le [1-(m-1)|t|]^{-\frac{1}{m-1}}-1.\)
In the scalar case $u_0 = 0$ and $t \geq 0$, equality holds identically, demonstrating again the sharpness of the comparison principle.
\end{example}


\begin{thebibliography}{9}
\bibitem{comtet} L. Comtet, \emph{Advanced Combinatorics: The Art of Finite and Infinite Expansions}, D. Reidel Publishing Company, Dordrecht, 1974.
\bibitem{gelfand} I. Gelfand, D. Raikov, and G. Shilov, \emph{Commutative Normed Rings}, Chelsea Publishing Company, New York, 1964.
\bibitem{orozco1} R. Orozco, \emph{Soluci\'on anal\'itica de la ecuaci\'on diferencial aut\'onoma de orden $k$}, Revista de Matem\'atica: Teor\'ia y Aplicaciones \textbf{23} (2016), no. 1, 63--84.
\bibitem{orozco2} R. Orozco, \emph{Ring of flows of one-dimensional differential equations
,} J. Algebra Appl. \textbf{21} 12 (2022), 2350004. https://doi.org/10.1142/S0219498823500044.
\end{thebibliography}
\end{document}